\documentclass[a4paper,oneside,12pt]{amsart}
\usepackage{geometry}               % See geometry.pdf to learn the layout options. There are lots.
\usepackage[english]{babel}
\usepackage[latin1]{inputenc}
\usepackage{amssymb,amsmath,amsthm,amsfonts}
\usepackage{graphicx}
\usepackage[dvipsnames]{xcolor}
\usepackage[colorlinks]{hyperref}
\usepackage{filecontents}
\usepackage[normalem]{ulem}

\numberwithin{equation}{section}

\allowdisplaybreaks

\newcommand{\N}{\mathbb{N}}
\newcommand{\Z}{\mathbb{Z}}
\newcommand{\R}{\mathbb{R}}

\newtheorem{theorem}{Theorem}[section]

\newtheorem{remark}[theorem]{Remark}
\theoremstyle{definition}

\renewcommand{\le}{\leqslant}

\renewcommand{\ge}{\geqslant}

\renewcommand{\epsilon}{\varepsilon}
\newcommand{\eps}{\varepsilon}
 \title[Fractional Burgers' equations]{Singularity formation\\ in fractional Burgers' equations}

\author{G. M. Coclite}\author{S. Dipierro}\author{F. Maddalena}\author{E. Valdinoci}
\address[Giuseppe Maria Coclite and Francesco Maddalena]{\newline
Department of Mechanics, Mathematics and Management, Polytechnic of Bari,
Via E.~Orabona 4, 70125 Bari, Italy.}
\email[]{giuseppemaria.coclite@poliba.it, francesco.maddalena@poliba.it}

\address[Serena Dipierro and Enrico Valdinoci]{\newline
Department of Mathematics and Statistics,
University of Western Australia,
35 Stirling Highway,
Crawley, WA 6009, Australia.}
\email[]{serena.dipierro@uwa.edu.au, enrico.valdinoci@uwa.edu.au}

\date{\today}


\thanks{The  authors are members of the Gruppo Nazionale per 
l'Analisi Matematica, la Probabilit\`a e le loro Applicazioni (GNAMPA) 
of the Istituto Nazionale di Alta Matematica (INdAM).
SD and EV are supported by the 
Australian Research Council
Discovery Project grant ``Nonlocal Equations at Work'' (NEW).
SV is supported by the DECRA Project ``Partial differential equations, free boundaries and applications''.}
\begin{document}
\maketitle

\begin{abstract}
The formation of singularities in finite time   in non-local Burgers' equations, with time-fractional derivative, is studied in detail. 
The occurrence of finite time singularity is proved, revealing the underlying mechanism, and precise estimates on the blow-up time are provided. 
The employment  of the present equation to model  a problem arising in job market is also
analyzed. 
\end{abstract}

\section{Introduction}

The study of singularities occurrence in nonlinear evolution problems constitutes a source of intriguing questions 
deeply related to the  mathematical and physical issues. The basic
 example of a PDE evolution leading to shock formation is given by the so called Burgers' equation (actually introduced by Airy \cite{A}) 
which represents a simple
model for studying the interaction between nonlinear and dissipative phenomena. 
Moreover, this equation exhibits the basic nonlinear mechanism shared by the more involved nonlinearities inherent to
Euler and Navier-Stokes equations~\cite{KS}. In exploiting a possible scenario for singularity formation in nonlocal 
evolution problems, continuing a line of research pursued in~\cite{CDMV}
from a different perspective, here we investigate the effect of a nonlocal
in time modification of Burgers' equation with respect to singularity creation. 
\medskip

Besides their interest from the purely mathematical point of view,
the nonlocal operators with respect to the time variable
find a number of concrete applications
in many emerging fields of research like, for instance,  the anomalous transportation problems
(see~\cite{BJW}),
the heat flow through ramified media (see~\cite{COMB}),
and the theory of viscoelastic fluids (see Section 10.2 in~\cite{MR1658022}
and the references therein).
Specifically, we will also present here
a concrete model from {\em job market} analysis
which naturally leads to a fractional Burgers' equation.
See also Chapter~1 in~\cite{CARB}
for several explicit motivations for fractional derivative problems.
\medskip

Focusing on the case of inviscid fluid mechanics, we recall that
in the classical Burgers' equation
explicit examples show the possible {\em formation of singularities
in finite time}, see \cite{B}. In particular, an initial condition with {\em unitary slope}
leads to a singularity at a {\em unitary time}.\medskip

The goal of this paper is to study whether a similar phenomenon
persists in nonlocal Burgers' equations with a time-fractional derivative.
That is, we investigate {\em how a memory effect
in the equation affects the singularity formation}.

Our main results are the following:
\begin{itemize}
\item The memory effect {\em does not prevent singularity formations}.
\item For initial data with unitary slopes, the blow-up time
can be explicitly estimated from {\em above}, in a way that is
uniform with respect to the memory effect (namely, it is {\em not
possible to slow down indefinitely the singularity formation
using only memory effects}).
\item Explicit bounds from {\em below} of the blow-up times are also possible.
\end{itemize}
\bigskip

The precise mathematical setting in which
we work is the following. First of all, to describe
memory effects,
we make use
of the left-Caputo-derivative of order~$\alpha\in(0,1)$
with initial time~$t_0$ for~$t\in(t_0,+\infty)$, defined by 
\begin{equation}\label{CAP} {}^C \! D^\alpha_{t_0,+} f(t):=
\frac{1}{\Gamma(1-\alpha)}
\int_{t_0}^t \frac{\dot f(\tau)}{(t-\tau)^\alpha}\,d\tau,\end{equation}
where~$\Gamma$ is the Euler Gamma function.

In this framework,
we consider the time-fractional Burgers' equation
driven by the left-Caputo-derivative, given by
\begin{equation}\label{BG} \begin{cases}
{}^C \! D^\alpha_{0,+} u(x,t) + u(x,t)\,\partial_x u(x,t)=0 & {\mbox{ for all $x\in\R$ and $t\in(0,T_\star)$,}}\\
u(x,0)=u_0(x).
\end{cases} \end{equation}
In the recent literature, various
types of fractional versions of the classical Burgers' equation
were taken into account from different
perspectives, see e.g.~\cite{MR1909206, MR2422665, MR3154615, MR3039000, MR3513173, MR3843183, MR3673503}
and the references therein (in this paper, we also
propose a simple motivation for
equation~\eqref{BG}
in Section~\ref{MOTIVAZ}).
When~$\alpha=1$, equation~\eqref{BG} reduces to the classical inviscid
Burgers' equation
\begin{equation}\label{SP}
\partial_t u(x,t) + u(x,t)\,\partial_x u(x,t)=0 .\end{equation}

\begin{remark}
{\rm We notice that examples of solutions to classical Burgers' equation, exhibiting instantaneous
 and spontaneous formation of singularities, work well also in the present case. Indeed, the aim of our study relies in understanding, through quantitative estimates, how the  
finite-time creation of singularities  can be affected by the presence of a fractional in time derivative.}
\end{remark}

We prove that the time-fractional Burgers' equation
driven by the left-Caputo-derivative may develop singularities
in finite time, according to the following result:

\begin{theorem}\label{CON} There exist a time~$T_\star>0$,
a function~$u_0\in C^\infty(\R)$,
a smooth solution~$u:\R\times[0,T_\star)\to\R$ of the time-fractional Burgers' equation in~\eqref{BG}
and a sequence~$t_n\nearrow T_\star$ as~$n\to+\infty$ such that
$$ \lim_{n\to+\infty} u(x,t_n)=\begin{cases}
+\infty & {\mbox{ if }} x\in(-\infty,0),\\
-\infty & {\mbox{ if }} x\in(0,+\infty).\end{cases}$$
\end{theorem}

\begin{remark}\label{R1} {\rm 
The function $u$ in Theorem~\ref{CON} will be constructed
by taking
\begin{equation}\label{MO} u(x,t):=-x\,v(t),\end{equation}
where~$v$ is the solution of the time-fractional equation
\begin{equation}\label{xvv-def}
\begin{cases}
{}^C \! D^\alpha_{0,+} v(t) = v^2(t) & {\mbox{ for }}t\in(0,T_\star),\\
v(0)=1.
\end{cases}
\end{equation}
When~$\alpha=1$, the equation in~\eqref{xvv-def}
reduces to~$\dot v(t)=v^2(t)$, which has the explicit
solution~$v(t)=1/(1-t)$. Therefore, the function in~\eqref{MO}
recovers the explicit, singular solution
\begin{equation}\label{CLA} u(x,t)=-\frac{x}{1-t}\end{equation}
of the classical Burgers' equation~\eqref{SP}
as~$\alpha\nearrow1$. Of course,
in the classical case, the blow-up time~$T_\star$
of~\eqref{CLA}
is exactly~$1$: in this sense,
our fractional construction in Theorem~\ref{CON},
recovers the classical case in the limit~$\alpha\nearrow1$.
}\end{remark}

We also observe that it is possible to give an explicit upper bound
on the blow-up time for the fractional solution~\eqref{MO}
in Theorem~\ref{CON}, as detailed
in the following result:

\begin{theorem}\label{ALTO}
If~$T_\star$ is the blow-up time found in Theorem~\ref{CON}, we have that
\begin{equation}\label{Bou1}
T_\star\le \left( \frac1{\Gamma(2-\alpha)}\right)^{1/\alpha}.
\end{equation}
In particular, for all~$\alpha\in(0,1)$,
\begin{equation}\label{Bou2}
T_\star\le e^{1-\gamma}=1.52620511\dots,
\end{equation}
where $\gamma$ is the Euler-Mascheroni constant.
\end{theorem}

\begin{remark}{\rm 
One can compare the general estimate in~\eqref{Bou2},
valid for all~$\alpha\in(0,1)$, with the blow-up time
for the classical solution in~\eqref{CLA}, in which~$T_\star=1$.
Indeed, we point out that the right hand side of~\eqref{Bou1}
approaches~$1$ as~$\alpha\nearrow1$. Hence,
in view of Remark~\ref{R1},
we have that the bound in~\eqref{Bou1} is optimal when~$\alpha\nearrow1$.
}\end{remark}

It is also possible to obtain a lower bound on the blow-up time
involving the right hand side in~\eqref{Bou1},
up to a reminder which is arbitrarily small as~$\alpha\nearrow1$.
Indeed, we have the following result:

\begin{theorem}\label{BASSO}
If~$T_\star$ is the blow-up time found in Theorem~\ref{CON}, we have that
for any~$\delta>0$ there exists~$c_\delta>0$ such that
\begin{equation}\label{Basso1}
T_\star\ge \frac{c_\delta^{\frac{1-\alpha}\alpha}}{1+\delta}
\left( \frac1{\Gamma(2-\alpha)}\right)^{1/\alpha}.
\end{equation}
\end{theorem}

\begin{remark} {\rm We observe that the right hand side of~\eqref{Basso1}
approaches~$1/(1+\delta)$ as~$\alpha\nearrow1$,
which, for small~$\delta$, recovers the unitary blow-up time
of the classical solution in~\eqref{CLA}.}
\end{remark}

\begin{remark} {\rm Of course, the blow-up time estimates
in Theorems~\ref{ALTO} and~\ref{BASSO} are specific for
the singular solution in~\eqref{MO}, and other singular solutions
have in general different blow-up times. As a matter of fact,
by scaling, if~$u$ is a solution of~\eqref{BG},
then so is~$u^{(\lambda)}(x,t):=u(\lambda^\alpha x,\lambda t)$,
for all~$\lambda>0$, with
initial datum~$u^{(\lambda)}_0(x):=u_0(\lambda^\alpha x)$.
In particular, if~$u$ is the function in~\eqref{MO}
and~$T_\star$ is its blow-up time,
then the blow-up time of~$u^{(\lambda)}$ is~$T_\star/\lambda$.
That is, when~$\lambda\in(1,+\infty)$, the slope of the initial datum
increases and accordingly the blow-up time becomes smaller. This is
the reason for which we choose
the setting in~\eqref{MO} to normalize the slope of the initial datum
to be unitary.}
\end{remark}

\begin{remark} {\rm It is interesting to observe
the specific effect of the Caputo derivative on the
solutions in simple and explicit
examples. {F}rom our perspective,
though the Caputo derivative is commonly viewed as a ``memory'' effect, the system does distinguish between
a {\em short-term} memory effect, which enhances the role
of the forcing terms, and a {\em long-term}
memory effect, which is more keen to remember the past configurations.

To understand our point of view on this phenomenon,
one can consider, for~$\alpha\in(0,1)$, the solution~$u=u(t)$ of the linear 
equation
\begin{equation}\label{ESEMPIOa}
\begin{cases}
{}^C \! D^\alpha_{0,+} u(t)=\displaystyle\sum_{k=1}^N \delta_{p_k}(t),\\
u(0)=0,\end{cases}\end{equation}
where~$0<p_1<\dots<p_N$ and~$\delta_p$ is the
Dirac delta at the point~$p\in\R$.

When~$\alpha=1$, equation~\eqref{ESEMPIOa}
reduces to the ordinary differential equation
with impulsive forcing term given by
\begin{equation}\label{ESEMPIO1}
\begin{cases}
\dot{u}(t)=\displaystyle \sum_{k=1}^N \delta_{p_k}(t),\\
u(0)=0.\end{cases}\end{equation}
Up to negligible sets,
the solution of~\eqref{ESEMPIO1} is the step function
\begin{equation}\label{SOLUZIONE1}
u(t) =\sharp \{ k\in\{1,\dots, N\}
{\mbox{ s.t. }} p_k<t\}=
\sum_{{1\le k\le N}\atop{p_k<t}} 1.
\end{equation}
On the other hand,
equation~\eqref{ESEMPIOa} is a Volterra-type problem
whose explicit solution is given by
\begin{equation}\label{SOLUZIONEa}
u(t) =\frac{1}{\Gamma(\alpha)}
\sum_{{1\le k\le N}\atop{p_k<t}} (t-p_k)^{\alpha-1}.
\end{equation}
Notice that the solution in~\eqref{SOLUZIONEa}
recovers~\eqref{SOLUZIONE1} as~$\alpha\nearrow1$.
Nevertheless, the sharp geometric difference between the solutions in~\eqref{SOLUZIONE1} and~\eqref{SOLUZIONEa} is apparent (see Figure \ref{fig:1}). 

\begin{figure}[htbp]
  \centering
  \includegraphics[width=0.7\linewidth]{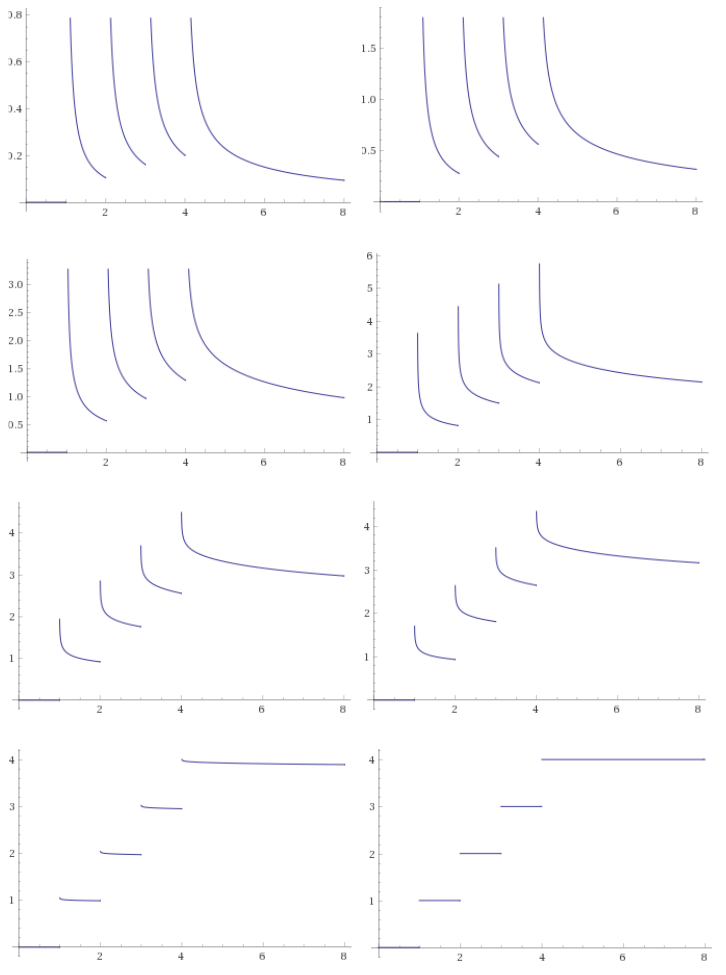}
  \caption{\it {{Plot of the solutions in~\eqref{SOLUZIONE1} and~\eqref{SOLUZIONEa}
    with the following parameters:
    $N=4$, $p_1=1$, $p_2=2$, $p_3=3$ and~$p_4=4$.
    The different plots correspond to the cases~$\alpha=\frac{1}{10}$, $\alpha=\frac{1}{4}$,
    $\alpha=\frac{1}{2}$, $\alpha=\frac{3}{4}$, $\alpha=\frac{7}{8}$, $\alpha=\frac{9}{10}$,
    $\alpha=\frac{99}{100}$, and~$\alpha=1$.}}}
    \label{fig:1}
\end{figure}

Indeed, while the classical solutions
experience a unit jump at the times where the impulses take place, the structure of
the fractional solutions exhibit a more complicated, and ``less monotone'', behavior. More specifically, on the one hand,
for fractional solutions, the short-term memory
effect of each impulse is to create a singularity towards
infinity, and in this sense its impact on the solution is much stronger than in the classical case. On the other hand, the solution in~\eqref{SOLUZIONEa} approaches zero outside the times in which the impulses occur, thus tending to recover the initial datum in view of a long-term memory effect.}
\end{remark}

The paper is organized as follows. 
Sections \ref{sec:2}, \ref{sec:3}, and \ref{sec:4} are devoted to the proofs of Theorems \ref{CON}, \ref{ALTO}, and \ref{BASSO}, respectively.
In Section \ref{MOTIVAZ} we propose a job market motivation for equation~\eqref{BG}.

\section{Proof of Theorem~\ref{CON}}
\label{sec:2}

The proof of Theorem~\ref{CON} relies
on a separation of variables method (as it will be apparent in the definition
of the solution~$u$ in~\eqref{u} at the end of this proof). To make
this method work, one needs a careful analysis of the solutions
of time-fractional equations, that we now discuss in details.
Fixed~$M\in\N\cap[ 4,+\infty)$, for any~$r\in\R$
we define~$f_M(r):=\min\{ r^2,M^2\}$. We let~$v_M$ be
the solution of the Cauchy problem
\begin{equation}\label{CP}
\begin{cases}
{}^C \! D^\alpha_{0,+} v_M(t) = f_M(v_M(t)) & {\mbox{ for }}t\in(0,+\infty),\\
v_M(0)=1.
\end{cases}
\end{equation}
The existence and uniqueness of the solution~$v_M$, which
is continuous up to~$t=0$,
is warranted by Theorem~2 on page~304 of~\cite{MR2252568}.
In addition, by Theorem~1 on page~300 of~\cite{MR2252568},
we know that this solution can be represented in an integral form by
the relation
\begin{equation}\label{VOLT}
v_M(t)=1+\frac1{\Gamma(\alpha)}\int_0^t \frac{f_M(v_M(\tau))}{(t-\tau)^{1-\alpha}}\,d\tau.
\end{equation}
In particular, since~$f_M\ge0$, we have that~$v_M\ge1$.
Also, by continuity at~$t=0$, there exists~$\delta>0$
such that
\begin{equation}\label{v0M}
{\mbox{$v_4(t)\le2$ for all~$t\in(0,\delta)$.}}\end{equation}
We claim that
\begin{equation}\label{v1M}
{\mbox{$v_M(t)=v_4(t)$ for all~$t\in(0,\delta)$ and all~$M\ge4$.}}
\end{equation}
Indeed, if~$t\in(0,\delta)$ and~$M\ge4$, we have that
$$ f_M(v_4(t))=\min\{ v_4^2(t),M^2\}=v_4^2(t)=\min\{ v_4^2(t),4^2\}=f_4(v_4(t)),$$
thanks to~\eqref{v0M}, and therefore~${}^C \! D^\alpha_{0,+} v_4(t) = f_M(v_4(t))$
for all~$t\in(0,\delta)$. Then, the uniqueness of the solution
of the Cauchy problem in~\eqref{CP} gives~\eqref{v1M}, as desired.

Furthermore, we observe that if~$M_2\ge M_1$ then~$f_{M_2}\ge f_{M_1}$
and then
$$ {}^C \! D^\alpha_{0,+} v_{M_2}(t) = f_{M_2}(v_{M_2}(t))
\ge f_{M_1}(v_{M_2}(t)).$$
Consequently, by the Comparison Principle\footnote{We observe
that we cannot use here the  Comparison Principle in Lemma~2.6
and Remark~2.1 on pages~219-220 in~\cite{MR3296607},
since the monotonicity of the nonlinearity goes in the opposite direction.}
in Theorem~4.10 on page~2894 in~\cite{MR3809535},
we conclude that~$v_{M_2}\ge v_{M_1}$.
Therefore, for every~$t\ge0$, we can define
\begin{equation}\label{v-1} v(t):=\lim_{M\to+\infty} v_M(t)=\sup_{M\in\N\cap[4,+\infty)} v_M(t)\in
[1,+\infty)\cup\{+\infty\}.\end{equation}
By~\eqref{v1M}, we know that
\begin{equation}\label{8} v(t)=v_4(t)\le \sup_{[0,\delta]} v_4<+\infty\qquad{\mbox{ for all~$t\in(0,\delta)$,}}\end{equation}
and hence we can consider the largest~$T_\star\in(0,+\infty)\cup\{+\infty\}$
such that
\begin{equation}\label{9} \sup_{t\in[0,T_0]}v(t)<+\infty\qquad{\mbox{ for all~$T_0\in(0,T_\star)$.}}\end{equation}
By~\eqref{8}, we have that~$T_\star\ge\delta$.
We claim that
\begin{equation}\label{v-def}
\begin{cases}
{}^C \! D^\alpha_{0,+} v(t) = v^2(t) & {\mbox{ for }}t\in(0,T_\star),\\
v(0)=1.
\end{cases}
\end{equation}
To prove this, we let~$T_0\in(0,T_\star)$ and we exploit~\eqref{9}
to see that
$$ M_0:=\sup_{t\in[0,T_0]}v(t)<+\infty,$$
and hence,
for every~$t\in(0,T_0)$ and every~$M\ge M_0$,
\begin{equation*}
f_M(v_{M_0}(t))=\min\{v_{M_0}^2(t), M^2\}
=v_{M_0}^2(t)=\min\{v_{M_0}^2(t), M_0^2\}=f_{M_0}(v_{M_0}(t)).
\end{equation*}
This gives that~${}^C \! D^\alpha_{0,+} v_{M_0}(t)=f_{M_0}(v_{M_0}(t))
=f_M(v_{M_0}(t))$ for all~$t\in(0,T_0)$ and~$M\ge M_0$, and therefore,
by the uniqueness of the solution
of the Cauchy problem in~\eqref{CP}, we find that~$v_M=v_{M_0}$
in~$(0,T_0)$. This and~\eqref{v-1} give that
$$ M_0\ge v(t)=v_{M_0}(t) \qquad{\mbox{for all }}t\in(0,T_0).$$
As a consequence, recalling~\eqref{VOLT},
we obtain that, for all~$t\in(0,T_0)$, the function~$v$
satisfies the integral relation
\begin{align*}
v(t)=&
v_{M_0}(t)=1+\frac1{\Gamma(\alpha)}\int_0^t \frac{f_{M_0}(v_{M_0}(\tau))}{
(t-\tau)^{1-\alpha}}\,d\tau\\
=&1+\frac1{\Gamma(\alpha)}\int_0^t \frac{f_{M_0}(v(\tau))}{
(t-\tau)^{1-\alpha}}\,d\tau=
1+\frac1{\Gamma(\alpha)}\int_0^t \frac{v^2(\tau)}{
(t-\tau)^{1-\alpha}}\,d\tau,
\end{align*}
and thus, by Theorem~1 in~\cite{MR2252568},
we obtain~\eqref{v-def}, as desired.

Now we claim that
\begin{equation}\label{TSI}
T_\star<+\infty.\end{equation}
To this end, we argue by contradiction and assume that~$T_\star=+\infty$.
We
let~$\lambda\ge2$, $T>0$ (which will be taken as large
as we wish in what follows), and
$$ \phi(t):=\begin{cases}
\left( 1-\displaystyle\frac{t}{T}\right)^\lambda & {\mbox{ if }}t\in[0,T],\\
0& {\mbox{ if }}t\in(T,+\infty).
\end{cases}$$
We know (see Lemmata~1 and~2 in~\cite{MR2464541}) that
\begin{equation}\label{L12}
\begin{split}
& \int^T_0 D^\alpha_{T,-}\phi(t)\,dt=\frac{\lambda\Gamma(\lambda-\alpha)}{
(\lambda-\alpha+1)\,\Gamma(\lambda-2\alpha+1)}\,T^{1-\alpha},
\\ &
\int^T_0 \frac{| D^\alpha_{T,-}\phi(t)|^2}{ \phi(t) }\,dt=\frac{\lambda^2}{\lambda+1-2\alpha}
\left( \frac{\Gamma(\lambda-\alpha)}{\Gamma(\lambda+1-2\alpha)}\right)^2
T^{1-2\alpha}.
\end{split}
\end{equation}
We also recall the left-Riemann-Liouville-derivative of order~$\alpha\in(0,1)$
with initial time~$t_0$ for~$t\in(t_0,+\infty)$, given by
$$ D^\alpha_{t_0,+} f(t):=
\frac{1}{\Gamma(1-\alpha)}\frac{d}{dt}
\int_{t_0}^t \frac{f(\tau)}{(t-\tau)^\alpha}\,d\tau,$$
and we point out that
\begin{equation*}
{}^C \! D^\alpha_{t_0,+} f(t)= D^\alpha_{t_0,+} \big( f(t)-f(t_0)\big).\end{equation*}
This and~\eqref{v-def} give that
\begin{equation}\label{C7}
v^2(t)=
{}^C \! D^\alpha_{0,+} v(t) =
D^\alpha_{0,+} \big( v(t)-v(0)\big)=D^\alpha_{0,+} w(t),
\end{equation}
where~$w(t):=v(t)-1$.

It is also useful to consider the
right-Riemann-Liouville-derivative of order~$\alpha\in(0,1)$
with final time~$t_0$ for~$t\in(-\infty,t_0)$, given by
\begin{equation*} D^\alpha_{t_0,-} f(t):=
-\frac{1}{\Gamma(1-\alpha)}\frac{d}{dt}
\int^{t_0}_t \frac{f(\tau)}{(\tau-t)^\alpha}\,d\tau.\end{equation*}
Integrating by parts
(see Corollary 2 on page~46 of~\cite{MR1347689},
or formula~(15) in~\cite{MR2728548}), and recalling~\eqref{C7},
we obtain that
\begin{equation}\label{PXuq}
\begin{split}
\int_0^T \phi(t)\,v^2(t)\,dt=&
\int_0^T \phi(t)\,D^\alpha_{0,+} w(t)\,dt\\
=&\int_0^T D^\alpha_{T,-} \phi(t)\,w(t)\,dt=
\int_0^T D^\alpha_{T,-} \phi(t)\,\big(v(t)-1\big)\,dt.\end{split}\end{equation}
{F}rom this and~\eqref{L12} we find that
\begin{equation*}
\int_0^T \phi(t)\,v^2(t)\,dt=
\int_0^T D^\alpha_{T,-} \phi(t)\,v(t)\,dt-C_1\,T^{1-\alpha},\end{equation*}
for some~$C_1>0$ independent of~$T$.

Furthermore,
\begin{align*}
\int_0^T D^\alpha_{T,-} \phi(t)\,v(t)\,dt=&
\int_0^T \frac{D^\alpha_{T,-} \phi(t)}{\sqrt{\phi(t)}}\,{\sqrt{\phi(t)}}\,v(t)\,dt\\
\le& \frac12\int_0^T \frac{|D^\alpha_{T,-} \phi(t)|^2}{{\phi(t)}}\,dt+
\frac12\int_0^T {{\phi(t)}}\,v^2(t)\,dt\\
=& C_2\,T^{1-2\alpha}+
\frac12\int_0^T {{\phi(t)}}\,v^2(t)\,dt,
\end{align*}
thanks to~\eqref{L12}, for some~$C_2>0$ independent of~$T$. As a consequence,
recalling~\eqref{PXuq},
we conclude that
\begin{equation*}
\frac12\int_0^T {{\phi(t)}}\,v^2(t)\,dt\le C_2\,T^{1-2\alpha}-C_1\,T^{1-\alpha}.
\end{equation*}
Therefore, recalling that~$v\ge1$ in view of~\eqref{v-1},
\begin{eqnarray*}&&
C_2\,T^{1-2\alpha}-C_1\,T^{1-\alpha} \ge
\frac12\int_0^T {{\phi(t)}}\,dt=
\frac12\int_0^T\left( 1-\displaystyle\frac{t}{T}\right)^\lambda\,dt=\frac{T}{2(1+\lambda)},
\end{eqnarray*}
and accordingly
\begin{eqnarray*}
0=\lim_{T\to+\infty} C_2\,T^{-2\alpha}-C_1\,T^{-\alpha}\ge\frac{1}{2(1+\lambda)},
\end{eqnarray*}
which is a contradiction, thus completing the proof of \eqref{TSI}.

Then, from~\eqref{9} and~\eqref{TSI}, we obtain that
\begin{equation*}
\limsup_{t\nearrow T_\star} v(t)=+\infty.
\end{equation*}
Hence, we consider a sequence~$t_n\nearrow T_\star$
such that
\begin{equation}\label{VB}
\lim_{n\to+\infty} v(t_n)=+\infty,
\end{equation}
and we define
\begin{equation}\label{u}
u(x,t):= -x\,v(t).
\end{equation}
For every~$t\in(0,T_\star)$, we have that
$$ {}^C \! D^\alpha_{0,+} u(x,t) + u(x,t)\,\partial_x u(x,t)=
-x\,{}^C \! D^\alpha_{0,+} v(t)+x\,v^2(t)=0,$$
thanks to~\eqref{v-def}, and also~$u(x,0)=-x\,v(0)=-x$.
These observations and~\eqref{VB}
prove Theorem~\ref{CON}.

\section{Proof of Theorem~\ref{ALTO}}
\label{sec:3}

We set
\begin{equation}\label{bal}
b=b(\alpha):=\left( \frac1{\Gamma(2-\alpha)}\right)^{1/\alpha}
\end{equation}
For any~$t\in(0,b)$, let also
$$ w(t):=\frac{b}{b-t}.$$
Notice that~$w(0)=1$. Moreover, for any~$t\in(0,b)$ and any~$\tau\in(0,t)$,
we have that
$$ \dot w(\tau)=\frac{b}{(b-\tau)^2} \le\frac{b}{(b-t)^2} =\frac{ w^2(t) }{b}
.$$
Consequently, by~\eqref{CAP}, for all~$t\in(0,b)$,
\begin{align*} {}^C \! D^\alpha_{0,+} w(t)=&
\frac{1}{\Gamma(1-\alpha)}
\int_{0}^t \frac{\dot w(\tau)}{(t-\tau)^\alpha}\,d\tau
\le\frac{w^2(t)}{b\,\Gamma(1-\alpha)}
\int_{0}^t \frac{d\tau}{(t-\tau)^\alpha}\\ =&
\frac{t^{1-\alpha} \,w^2(t)}{b\,\Gamma(1-\alpha)\,(1-\alpha)}=
\frac{t^{1-\alpha} \,w^2(t)}{b\,\Gamma(2-\alpha)}\le
\frac{b^{1-\alpha} \,w^2(t)}{b\,\Gamma(2-\alpha)}\\=&
\frac{w^2(t)}{b^\alpha\,\Gamma(2-\alpha)}=w^2(t).
\end{align*}
Therefore, using the
Comparison Principle in Theorem~4.10 on page~2894 in~\cite{MR3809535},
if~$v$ is as in~\eqref{v-def}, we find that~$v\ge w$
in their common domain of definition. This, \eqref{u},
and the fact that~$w$ diverges at~$t=b$ yield that
\begin{equation}\label{Bou3}
T_\star\le b=b(\alpha),\end{equation}
which, together with~\eqref{bal}, establishes~\eqref{Bou1}, as desired.

Now we prove~\eqref{Bou2}. For this, we first show that
the map~$(0,1)\ni\alpha\mapsto b(\alpha)$ that was introduced in~\eqref{bal}
is monotone.
To this end, we recall the polygamma functions for~$\tau\in(1,2)$
and~$n\in\N$ with their integral representations,
namely
\begin{align*}
&\psi _{n}(\tau):=
\left({\frac {d}{d\tau}}\right)^{n+1}
\log (\Gamma (\tau))\qquad{\mbox{ for all }}n\in\{0,1,2,3,\dots\},\\
&(-1)^{n+1}\int _{0}^{\infty }{\frac {t^{n}e^{-t\tau}}{1-e^{-t}}}\,dt
\qquad{\mbox{ for all }}n\in\{1,2,3,\dots\}.\end{align*}
We observe, in particular, that, for all~$\tau\in(1,2)$,
\begin{equation}\label{psi1}
\psi_1(\tau)=\int _{0}^{\infty }{\frac {t e^{-t\tau}}{1-e^{-t}}}\,dt>0.
\end{equation}
Let also, for all~$\tau\in(1,2)$,
$$ \xi(\tau):=\log(\Gamma(\tau))+(2-\tau)\psi_0(\tau).$$
We see that
$$ \xi'(\tau)=\psi_0(\tau)-\psi_0(\tau)+(2-\tau)\psi_1(\tau)>0,$$
thanks to~\eqref{psi1} and therefore, for all~$\tau\in(1,2)$,
\begin{equation}\label{xi}
0<\int_\tau^2\xi'(\sigma)\,d\sigma=\xi(2)-\xi(\tau)=
\log(\Gamma(2))-\xi(\tau)=-\xi(\tau).\end{equation}
Now we define
$$ \lambda(\tau):=\frac{\log(\Gamma(\tau))}{2-\tau}.$$
We have that, for all~$\tau\in(1,2)$,
$$ \lambda'(\tau)=\frac{\log(\Gamma(\tau))}{(2-\tau)^2}+
\frac{\psi_0(\tau)}{2-\tau}=\frac{\xi(\tau)}{(2-\tau)^2}<0,$$
due to~\eqref{xi}. 

Therefore, the function~$(1,2)\ni\tau\mapsto\lambda(\tau)$ is
decreasing, and hence so is the function~$(1,2)\ni\tau\mapsto
e^{\lambda(\tau)}=:\Lambda(\tau)$. Hence,
using the substitution~$\tau:=2-\alpha$,
with~$\alpha\in(0,1)$, we deduce that
the following function is increasing:
\begin{align*} \Lambda(2-\alpha)=&
e^{\lambda(2-\alpha)}=
\exp\left(\frac{\log(\Gamma(2-\alpha))}{\alpha}\right)\\=&
\exp\left({\log (\Gamma^{1/\alpha}(2-\alpha))}\right)=
\Gamma^{1/\alpha}(2-\alpha)=\frac{1}{b(\alpha)},
\end{align*}
thanks to~\eqref{bal}.

Consequently, the function~$(0,1)\ni\alpha\mapsto b(\alpha)$
is decreasing, hence it attains its maximum as~$\alpha\searrow0$.
This and~\eqref{Bou3} give that
\begin{equation}\label{Bou4}
T_\star\le \lim_{\alpha\searrow0} b(\alpha).\end{equation}
Furthermore, using L'H\^{o}pital's Rule,
\begin{eqnarray*}&&
\lim_{\alpha\searrow0}\frac{\log(\Gamma(2-\alpha))}{\alpha}=
-\lim_{\alpha\searrow0}\psi_0(2-\alpha)=-\psi_0(2)=\gamma-1,
\end{eqnarray*}
where~$\gamma$ is the Euler-Mascheroni constant, and therefore
$$ \lim_{\alpha\searrow0}b(\alpha)=
\lim_{\alpha\searrow0}
\exp\left(-\frac{\log(\Gamma(2-\alpha))}{\alpha}\right)
=e^{1-\gamma}.$$
This and~\eqref{Bou4} give the desired result in~\eqref{Bou2}.

\section{Proof of Theorem~\ref{BASSO}}
\label{sec:4}

We let~$\delta>0$ as in the statement of Theorem~\ref{BASSO}, and
\begin{equation}\label{kappa}
\kappa:=\sqrt{1+\delta}-1>0.\end{equation}
We define
\begin{align*}
\eta:=\frac{(1+\kappa)^2}{\kappa^2},&\qquad
d:=\left( \frac{1}{\Gamma(2-\alpha)\,\kappa\,\eta\,(1+\eta)}\right)^{\frac1\alpha},\\
 a:=\frac{\Gamma(2-\alpha)}{d^{1-\alpha}},&
%% =\big(\Gamma(2-\alpha)\big)^{\frac1\alpha}
%% \big(\kappa\,\eta\,(1+\eta)\big)^{\frac{1-\alpha}\alpha}
\qquad
b:=(1+\kappa)a.
\end{align*}
Let also
\begin{equation}\label{TI} T:= \frac1b-(1+\eta)d.\end{equation}
In light of~\eqref{TI}, we remark that
\begin{equation}\label{QQ}
\begin{split}
T\,&=\frac{1}{(1+\kappa)a}-(1+\eta)d\\
&=\frac{d^{1-\alpha}}{(1+\kappa)\,\Gamma(2-\alpha)}-(1+\eta)d\\
&=(1+\eta)\,d\,\left(
\frac{d^{-\alpha}}{(1+\kappa)\,\Gamma(2-\alpha)\,(1+\eta)}-1
\right)\\
&=(1+\eta)\,d\,\left(
\frac{\kappa\eta}{1+\kappa}-1\right)\\
&=(1+\eta)\,d\,\left(
\frac{1+\kappa}{\kappa}-1\right)\\
&=\frac{(1+\eta)\,d}{\kappa}\\
&=\frac{1+\eta}{\kappa}\;
\left( \frac{1}{\Gamma(2-\alpha)\,\kappa\,\eta\,(1+\eta)}\right)^{\frac1\alpha}
\\&=
\frac{1}{\big( \Gamma(2-\alpha)\big)^{\frac1\alpha}\,\kappa^{\frac{1+\alpha}\alpha}\,\eta^{\frac1\alpha}\,(1+\eta)^{\frac{1-\alpha}\alpha}}
\\&=
\frac{\kappa^{\frac{3(1-\alpha)}\alpha}}{
\big( \Gamma(2-\alpha)\big)^{\frac1\alpha}\,(1+\kappa)^{\frac2\alpha}\,(1+2\kappa+2\kappa^2)^{\frac{1-\alpha}\alpha}}
.\end{split}
\end{equation}
Recalling~\eqref{kappa}, we can also define
$$ c_\delta:=
\frac{\kappa^{3}}{(1+\kappa)^{2}\,(1+2\kappa+2\kappa^2)},
$$
and then~\eqref{QQ} becomes
\begin{equation}\label{QQ2}
T=\frac{c_\delta^{\frac{1-\alpha}\alpha}}{
\big( \Gamma(2-\alpha)\big)^{\frac1\alpha}\,(1+\kappa)^{2}}=
\frac{c_\delta^{\frac{1-\alpha}\alpha}}{
\big( \Gamma(2-\alpha)\big)^{\frac1\alpha}\,(1+\delta)},
\end{equation}
which coincides with the right hand side of~\eqref{Basso1}.

Therefore, to complete the proof of
Theorem~\ref{BASSO}, it is enough to show that
\begin{equation}\label{FG}
T_\star\ge T.
\end{equation}
To this end, for all~$t\in(0,T)$, we define
$$ z(t):=\frac{b}{a(1-bt)}+1-\frac{b}{a}.
$$
Notice that~$z(0)=1$. Moreover, for all~$t\in(0,T)$
and~$\tau\in(0,t)$ we have that
\begin{equation}\label{BAA}
\begin{split}
\frac{b}{a(1-b\tau)}-z(\tau+d)\,&=
\frac{b}{a(1-b\tau)}-\frac{b}{a(1-bd-b\tau)}-1+\frac{b}{a}\\
&=
\frac{1+\kappa}{ 1-b\tau }-\frac{1+\kappa}{1-bd-b\tau }+\kappa\\
&=
-\frac{(1+\kappa)\,b\,d}{ (1-b\tau)(1-bd-b\tau) }+\kappa\\
&\ge
-\frac{(1+\kappa)\,b\,d}{ (1-bT)(1-bd-bT) }+\kappa\\
&=
-\frac{(1+\kappa)}{ (1+\eta)\,\eta\,b\,d }+\kappa.
\end{split}\end{equation}
Hence, since
$$ bd=(1+\kappa)ad=
(1+\kappa)\,\Gamma(2-\alpha)\,d^{\alpha}=
\frac{(1+\kappa)}{\kappa\,\eta\,(1+\eta)},$$
we see from~\eqref{BAA} that
\begin{eqnarray*}
\frac{b}{a(1-b\tau)}-z(\tau+d)\ge
-\kappa+\kappa=0,
\end{eqnarray*}
and therefore
$$ \frac{b^2}{a^2(1-b\tau)^2}\ge z^2(\tau+d).$$
As a consequence, we conclude that
\begin{eqnarray*}
\dot z(\tau)=\frac{b^2}{a(1-b\tau)^2}\ge a\,z^2(\tau+d).
\end{eqnarray*}
Accordingly, by~\eqref{CAP}, for all~$t\in(0,T)$,
\begin{eqnarray*}
{}^C \! D^\alpha_{0,+} z(t)&=&
\frac{1}{\Gamma(1-\alpha)}
\int_{0}^t \frac{\dot z(\tau)}{(t-\tau)^\alpha}\,d\tau\\
&\ge&\frac{a}{\Gamma(1-\alpha)}
\int_{t-d}^t \frac{ z^2(\tau+d)}{(t-\tau)^\alpha}\,d\tau.
\end{eqnarray*}
Consequently, using the fact that~$z$ is increasing,
$$ {}^C \! D^\alpha_{0,+} z(t)\ge
\frac{a}{\Gamma(1-\alpha)}
\int_{t-d}^t \frac{ z^2(t)}{(t-\tau)^\alpha}\,d\tau=
\frac{a\,d^{1-\alpha}\,z^2(t)}{\Gamma(2-\alpha)}=z^2(t).
$$
Then, recalling~\eqref{v-def}
and exploiting the Comparison Principle
in Theorem~4.10 on page~2894 in~\cite{MR3809535},
we obtain that~$v\le z$ in their common domain of definition.
In particular, this gives~\eqref{FG},
and so the proof of Theorem~\ref{BASSO} is complete.

\section{A motivation for~(\ref{BG}) from the job market}
\label{MOTIVAZ}

In this section we give a simple, but concrete, motivation
for the time-fractional Burgers' equation
in~\eqref{BG} making a model of an ideal job market
from a few basic principles. The discussion that we present
here is a modification of classical models proposed
for fluid dynamics and traffic flow in a highway.

We fix parameters~$\delta$, $\eps>0$
and we use the real line to describe the positions
available in a company, in which workers can decide to work.
More specifically, the working levels in the company are
denoted by~$x\in\eps\Z$ and the higher the value of~$x$
the higher and more appealing the position is
(e.g., $x=\eps$ corresponds to Brigadier,
$x=2\eps$ to Major, $x=3\eps$
to Lieutenant, $x=4\eps$ to General, etc.).

We suppose that the main motivation for a worker to join the company
by taking the position~$x\in\eps\Z$ at time~$t\in\delta\N$
is provided by the possibility
of career progression towards the successive level.
If we denote by~$\rho$ the number of people
employed in a given position at a given time,
and by~$v$ the velocity of career progression
relative to a given position at a given time,
the ``group velocity'' of career progression for
a given position at a given time
is obtained by the product~$p:=\rho v$. 

We suppose that
the potential worker who is possibly entering
the company at the level~$x\in\eps\Z$
will look at the value of~$p$ for its perspective position
and compare it with the value of~$p$
relative to subsequent level~$x+\eps$, and this will constitute,
in this model, the main drive for the worker to join the company.
At time~$t\in\delta\N$, this  driving force is therefore quantified by
\begin{equation} \label{LADD}
{\mathcal{D}}(x,t):=
\hat{c}\big(p(x+\eps,t)-p(x,t)\big)=\hat{c}\Big(\rho(x+\eps,t)\,v(x+\eps,t)-
\rho(x,t)\,v(x,t)\Big), \end{equation}
for a normalizing constant~$\hat{c}>0$
Then, we assume that the potential worker
base her or his decision not only considering the driving force
at the present time, but also taking into
account the past history of the company. Past events
will be weighted by a kernel~${\mathcal{K}}$, to make the information coming from
remote times less important than the ones relative to the contemporary
situation. For concreteness, we suppose that
the information coming from the time~$t-\tau$, with~$t=\delta N$, $N\in\N$,
and~$\tau\in\{\delta,2\delta,\dots, \delta N\}$, is weighted by the kernel
\begin{equation}\label{KER}
{\mathcal{K}}(\tau):=\frac{\delta}{\tau^\beta},\qquad
{\mbox{ for some }}\beta\in(0,1).
\end{equation}
If all the potential workers argue in this way, the number of workers
at time~$t=\delta N$ in the working position~$x\in\eps\Z$
of the company is given by the initial number of workers,
incremented by the effect of the drive function
in the history of the company, according to the memory effect that
we have described, that is
\begin{equation*}
\rho(x,t)=\rho(x,\delta N)=\rho(x,0)+c\,
\sum_{j=1}^N {\mathcal{D}}(x,t-\delta j) \,{\mathcal{K}}(\delta j),
\end{equation*}
for some normalizing constant~$c>0$.
Hence, exploiting~\eqref{KER},
\begin{equation}\label{SUM}
\rho(x,t)=\rho(x,0)+c\,
\sum_{j=1}^N {\mathcal{D}}(x,t-\delta j) \,
\frac{\delta}{(\delta j)^\beta}.
\end{equation}
Using the Riemann sum approximation of an integral,
for small~$\delta$ we can substitute the summation in the right hand side
of~\eqref{SUM} with an integral, and, with this asymptotic
procedure, we replace~\eqref{SUM} with
\begin{equation}\label{IN}
\rho(x,t)=\rho(x,0)+c\,
\int_0^t {\mathcal{D}}(x,t-\tau) \,
\frac{d\tau}{\tau^\beta}.
\end{equation}
Then, we define~$\alpha:=1-\beta\in(0,1)$ and, up to a time scale,
we choose~$c:=1/\Gamma(\alpha)$. In this way, we can write~\eqref{IN}
as
\begin{align*}
\rho(x,t) =&\rho(x,0)+\frac1{\Gamma(\alpha)}\,
\int_0^t {\mathcal{D}}(x,t-\tau) \,
\frac{d\tau}{\tau^{1-\alpha}}\\
=&\rho(x,0)+\frac1{\Gamma(\alpha)}\,
\int_0^t {\mathcal{D}}(x,\sigma) \,
\frac{d\sigma}{(t-\sigma)^{1-\alpha}},
\end{align*}
or, equivalently (see e.g. Theorem~1 on page~300 of~\cite{MR2252568}),
\begin{equation*}
{}^C \! D^\alpha_{0,+} \rho(x,t)={\mathcal{D}}(x,t).
\end{equation*}
Thus, recalling~\eqref{LADD} (and using the normalization~$\hat{c}:=1/\eps$),
\begin{equation*}
{}^C \! D^\alpha_{0,+} \rho(x,t)=
\frac{\rho(x+\eps,t)\,v(x+\eps,t)-
\rho(x,t)\,v(x,t)}{\eps},\end{equation*}
and then, in the approximation of~$\eps$ small,
\begin{equation}\label{vcc1}
{}^C \! D^\alpha_{0,+} \rho(x,t)=\partial_x\Big(
\rho(x,t)\,v(x,t)\Big).\end{equation}
Now we make the ansatz that the career velocity is mainly influenced
by the number of people in a given position, namely this velocity
is proportional to the ``vacancies'' in a given working level.
If~$\rho_{\max}\in(0,+\infty)$ is the maximal number of workers 
that the market allows in any given position, we therefore assume that
\begin{equation}\label{vcc} v=\tilde{c} (\rho_{\max}-\rho),\end{equation}
for a normalizing constant~$\tilde{c}>0$. Of course, in more
complicated models, one can allow~$\rho_{\max}$ and~$\tilde{c}$
to vary in space and time, but we will take them to be constant
to address the simplest possible case, and in fact, for simplicity,
up to scalings,
we take~$\rho_{\max}=1$ and~$\tilde{c}=1$.

Then, plugging~\eqref{vcc} into~\eqref{vcc1}, we obtain
\begin{equation}\label{RHOBURG} {}^C \! D^\alpha_{0,+} \rho(x,t)= \partial_x\Big(
\rho(x,t)\,\big(1-\rho(x,t)\big)\Big).\end{equation}
Now we perform the substitution
\begin{equation}\label{RHO-U-BU} u(x,t):=2\rho(x,t)-1,\end{equation}
and we thereby conclude that
\begin{align*}
{}^C \! D^\alpha_{0,+}u(x,t)=&
2{}^C \! D^\alpha_{0,+} \rho(x,t)\\
=&2\partial_x\Big(
\rho(x,t)\,\big(1-\rho(x,t)\big)\Big)\\
=&2\partial_x\left(
\frac{u(x,t)+1}{2}\,\left(1-\frac{u(x,t)+1}{2}\right)\right)
\\=&2\partial_x\left(
\frac{u(x,t)+1}{2}-\frac{u^2(x,t)+2u(x,t)+1}{4}\right)
\\=&2\partial_x\left(
\frac{1}{4}-\frac{u^2(x,t)}{4}\right)\\
=&-u(x,t)\,\partial_x u(x,t),
\end{align*}
which corresponds to~\eqref{BG}.
\medskip

We remark that in the model above one can interpret~$\rho\in\R$ also when it takes negative values,
e.g. as a position vacancy.
As a matter of fact, since
the driving force of equation~\eqref{RHOBURG}
can be written as~$\partial_x\rho(1-2\rho)$,
we observe that such a drive becomes ``stronger''
for negative values of~$\rho$ (that is, vacancies
in the job markets tend to increase the number of filled positions).
\medskip

It is also interesting to interpret the result in Theorem~\ref{CON}
in the light of the motivation discussed here
and recalling the setting in~\eqref{RHO-U-BU}.
Indeed, the value~$1/2$ for the working force~$\rho$
plays a special role in our framework since not only it
corresponds to the average between the null working force and the maximal one allowed by the market, but also, and most importantly, to the critical value of the concave function~$\rho(1-\rho)$,
whose derivative is the driving force of equation~\eqref{RHOBURG}.

In this spirit, recalling~\eqref{MO},
we have that the solution found in Theorem~\ref{CON}
takes the form
\begin{equation}\label{rhjsmarke} \rho(x,t)=\frac{1-xv(t)}{2},\end{equation}
for a function~$v$ which is diverging in finite time.
The expression in~\eqref{rhjsmarke} says that the role corresponding to the job position~$x=0$ has, at the initial time, exactly the critical working force~$\rho=1/2$. Given the linear structure
in~$x$ of the solution in~\eqref{rhjsmarke}, this says that the job position corresponding to~$x=0$ will maintain its critical value~$\rho=1/2$ for all times, while higher level job roles will experience a dramatic loss of number of positions available (and, correspondingly, lower level job roles a dramatic increase). Though it is of course unrealistic that the job market really attains an (either positive or negative) infinite value in a finite time, and the model presented in equation~\eqref{RHOBURG} must necessarily ``break'' for too large values
of~$\rho$ (which, of course, in practice,
cannot exceed the total
working population),
we think that solutions
such as~\eqref{rhjsmarke} may represent a concrete case
in which the market would in principle allow
arbitrarily high level job positions, but in practice (almost)
all the workers end up obtaining a position level below a certain threshold (in this case normalized to~$x=0$),
which constitutes a``de facto" optimal role allowed by the
evolution of special preexisting conditions.

\bibliographystyle{abbrv}
\bibliography{CB-ref}

\vfill
\end{document}